\RequirePackage{etoolbox}
\csdef{input@path}{{style/}{graphics/}}
\documentclass[aos,MSNbibl,seceqn,dvips]{arximspdf}

\doi{10.1214/14-AOS1289}
\volume{43}
\issue{2}
\pubyear{2015}
\firstpage{652}
\lastpage{674}
\docsubty{FLA}

\makeatletter
\newcommand{\lleft}{\left}
\newcommand{\rrvert}{\vert}
\newcommand{\rright}{\right}
\newcommand{\rrVert}{\Vert}
\newcommand{\llvert}{\vert}
\newcommand{\llVert}{\Vert}
\renewcommand{\mid}{|}
\newcommand{\implies}{\Longrightarrow}
\newcommand{\eqref}[1]{(\ref{#1})}
\newcommand{\bbR}{\mathbb{R}}
\newcommand{\bbN}{\mathbb{N}}
\newtheorem{theorem}{Theorem}[section]
\newtheorem{lemma}[theorem]{Lemma}
\newproclaim{remark}[theorem]{Remark}
\newproclaim{ass}{Assumption}
\newcommand{\gp}{\operatorname{GP}}
\newcommand{\bbH}{\mathbb{H}}
\newcommand{\bbB}{\mathbb{B}}
\makeatother

\begin{document}
\begin{frontmatter}

\title{Minimax-optimal nonparametric regression in high~dimensions}
\runtitle{Minimax-optimal nonparametric regression}

\begin{aug}
\author[A]{\fnms{Yun}~\snm{Yang}\corref{}\ead[label=e1]{yy84@berkeley.edu}}
\and
\author[B]{\fnms{Surya T.}~\snm{Tokdar}\ead[label=e2]{tokdar@stat.duke.edu}}
\runauthor{Y. Yang and S.~T. Tokdar}
\affiliation{University of California, Berkeley
and Duke University}
\address[A]{Department of EECS\\
University of California, Berkeley\\
Berkeley, California 94720\\
USA\\
\printead{e1}}
\address[B]{Department of Statistical Science\\
Duke University\\
Box 90251\\
Durham, North Carolina 27708-0251\\
USA\\
\printead{e2}}
\end{aug}

%
\received{\smonth{8} \syear{2014}}
%
\revised{\smonth{11} \syear{2014}}

%
\begin{abstract}
Minimax $L_2$ risks for high-dimensional nonparametric regression are
derived under two sparsity assumptions: (1) the true regression surface
is a sparse function that depends only on $d=O(\log n)$ important
predictors among a list of $p$ predictors, with $\log p= o(n)$; (2) the
true regression surface depends on $O(n)$ predictors but is an additive
function where each additive component is sparse but may contain two or
more interacting predictors and may have a smoothness level different
from other components. For either modeling assumption, a practicable
extension of the widely used Bayesian Gaussian process regression
method is shown to adaptively attain the optimal minimax rate (up to
$\log n$ terms) asymptotically as both $n,p \to\infty$ with $\log p = o(n)$.
\end{abstract}

%
\begin{keyword}[class=AMS]
\kwd[Primary ]{62G08}
\kwd{62C20}
\kwd[; secondary ]{60G15}
\end{keyword}
\begin{keyword}
\kwd{Adaptive estimation}
\kwd{high-dimensional regression}
\kwd{minimax risk}
\kwd{model selection}
\kwd{nonparametric regression}
\end{keyword}
\end{frontmatter}

\section{Introduction}\label{sec1}

Rapid advances in technology have empowered\break researchers to collect data
on a large number of explanatory variables to predict many outcomes of
interest \cite{peter2011}. Because the relationship between an outcome
$Y$ and its predictors $X_1, \ldots, X_p$ may be highly nonlinear and
involve interaction, there is a practical need to investigate
statistical estimation under multivariate regression models
%
\begin{equation}
\label{eqreg} Y= \mu+ f(X_1,\ldots,X_p)+\varepsilon,\qquad
\varepsilon\sim N\bigl(0,\sigma^2\bigr),
\end{equation}
with minimal assumptions made on $f$. The quality of estimation that
may be achieved under an assumed model can be mathematically quantified
by the minimax risk of estimating $f$ from $n$ data points. A classic
result due to Charles Stone \cite{Stone1982} states that if no
assumption is made on how $f$ depends on $X_1,\ldots,X_p$ other than
requiring it to be differentiable with a smoothness level $\alpha> 0$
(definition below), then the associated minimax risk decays in $n$ at a
rate $n^{-\alpha/ (2\alpha+ p)}$. This rate is very slow when $p$ is
large, which means a very large sample size is needed for high quality
statistical estimation---a phenomenon that has been termed the ``curse
of dimensionality.'' The curse of dimensionality becomes even more
pronounced in the so-called large $p$ small $n$ setting, where the
minimax risk decay rate is expressed as a function of both $n$ and $p$,
with $p$ growing faster than $n$ \cite{Verzelen2012}.

Practically motivated modeling assumptions must focus on nonparametric
spaces of functions with lower inherent dimensions than the manifest
dimension~$p$. An example of such assumptions is
\begin{enumerate}[M1.]
\item[M1.] $f$ potentially depends on all elements of $X = (X_1,
\ldots, X_p)$, but $X$ itself lies in a low-dimensional manifold
$\mathcal{M}^d$ in the ambient space $\bbR^p$.
\end{enumerate}
It is well known that under M1, the minimax rate is $n^{-2\alpha/
(2\alpha+ d)}$ which is determined by the smoothness level $\alpha$
of $f$ and the latent manifold dimension $d$ \cite
{Bickel2007,Ye2008,Yang2013,Scott2006,Kulkarni1995,Kpotufe2012}, but
does not depend on
the ambient dimension $p$. Various nonparametric regression techniques
that operate on the ambient space and do not require estimation of the
underlying manifold indeed achieve this minimax rate without any prior
knowledge of $d$ or $\alpha$ \cite{Kulkarni1995,Yang2013}.

However, for many high-dimensional applications, such as gene
expression studies, a low-dimensional manifold assumption on $X$ may
not be tenable or verifiable. In such cases, one often uses the
following sparsity inducing assumption:
\begin{enumerate}[M2.]
\item[M2.] $f$ depends on a small subset of $d$ predictors with $d \le
\min\{n,p\}$.
\end{enumerate}
M2 has served as the springboard for many widely used regression
methods, including high-dimensional linear regression approaches, such
as the Lasso \cite{Tibshirani1996} and the Dantzig selector \cite
{candes2007}, and nonparametric regression methods with variable
selection, such as the Rodeo \cite{Lafferty2008} and the Gaussian
process regression \cite{Rasmussen2006}. The latter two allow flexible
shape estimation of $f$ and is able to capture interactions among the
selected important predictors. However, in light of the classic result
due to \cite{Stone1982} it is conceivable that when $f$ is allowed to
be fully nonparametric, M2~should also suffer from the curse of
dimensionality in a large $p$ small $n$ setting, unless $d$ is much
much smaller than $p$, that is, the regression function is assumed to
be extremely sparse. A precise result that extends the work of \cite
{Stone1982} to account for predictor selection is presented in
Section~\ref{sminimax}.

To relax this assumption of extreme sparsity without having to
completely give up on nonparametric shape flexibility, we introduce a
third modeling assumption:
\begin{enumerate}[M3.]
\item[M3.] $f$ may depend on $d\asymp\min\{n^\gamma,p\}$ variables
for some $\gamma\in(0,1)$ but admits an additive structure $f=\sum
_{s=1}^kf_s$, where each component function $f_s$ depends on a small
$d_s$ number of predictors.
\end{enumerate}
Clearly, M3 subsumes M2 as a special case. In Section~\ref{sminimax},
we show M2 gives slowest minimax rates within M3. At the opposite
extreme is the modeling assumption that $f$ admits a completely
additive structure with univariate components $f(X) = f_1(X_{i_1}) +
\cdots+ f_d(X_{i_d})$ for which scalable algorithms have been devised
\cite{Hastie1986} and attractive minimax risk bounds have been derived
albeit under the strong assumption that all component functions $f_s$
have the same smoothness level \cite
{Koltchinskii2010,Meier2009,Ravikumar2009,Raskutti2012}.

Compared to either of these two extremes, M3 provides a much more
practically attractive theory of large $p$ nonparametric regression. In
Theorem \ref{thmM3}, we derive sharp upper and lower bounds on the
minimax $L_2$ estimation risk under M3 as a function of $n$, $p$, $k$,
component sizes $d_1,\ldots,d_k$ and smoothness levels of $f_1,\ldots
,f_k$ which are allowed to have different levels of smoothness than one
another. Minimax rates under M2 and the completely additive structure
of \cite{Raskutti2012} follow as corollaries to this general result.
Our calculations suggest that M3 offers a minimax risk that decays
quickly in $n$ even when $p$ grows almost exponentially in $n$, $f$
involves nearly $\log p$ many predictors and these predictors interact
with each other.

In Section~\ref{sGP}, we demonstrate that a conceptually
straightforward extension of the widely used Gaussian process
regression method (see, e.g., \cite{Rasmussen2006}, for a
review)
achieves the minimax rate adaptively across all subclasses of M3 under
suitable large $p$ small $n$ asymptotics where $p$ grows almost
exponentially in $n$. In this paper we restrict only to a theoretical
study of this new approach, which we name ``additive Gaussian process
regression.'' This approach appears entirely practicable with
computational demands similar to those of the popular Bayesian additive
regression tree method \cite{bartpaper}. A full fledged
methodological development of the same is underway and will be reported
elsewhere.

The rest of the paper is organized as follows. Section~\ref{snotations}
introduces the notation and some basic assumptions.
Section~\ref{sminimax} summarizes our main minimax results for
high-dimensional
nonparametric regression under M2 and M3. Section~\ref{sGP} proves the
adaptive minimax optimality of additive Gaussian
process regression. Section~\ref{sproofs} provides proofs of our main
results in Sections~\ref{sminimax}~and~\ref{sGP}. Supporting
technical results and proofs are presented in Section~\ref{saux}.

\section{Notation}\label{sec2}
\label{snotations}
Let $(X^i, Y^i)$, $i = 1, \ldots, n$ denote the observations on $(X, Y)$.
We make a stochastic design assumption that $X^1, \ldots, X^n$ are
independent and identically distributed (IID) according to some
compactly supported probability measure $Q$ on $\bbR^p$ and that $f
\in L_2(Q)$, the linear space of real valued functions on $\bbR^p$
equipped with inner product $\langle f, g\rangle_Q = \int f(x)g(x)
Q(dx)$ and norm $\llVert f\rrVert _Q = \langle f, f\rangle_Q^{1/2}$.
We\vspace*{1pt} do not
need to know or estimate $Q$ for the purpose of estimating $f$, but it
is a natural candidate to judge average prediction accuracy at future
observations of $X$ drawn from $Q$, as will be the case under simple
exchangeability assumptions. Without loss of generality assume support~$(Q) \subset[0,1]^p$. Let $\llVert \cdot\rrVert $ stand for the $L_2$
norm under the Lebesgue measure.

The $L_2$ minimax risk of estimation associated with any function space
$\Sigma\subset L_2(Q)$ is defined as
\[
r^2_n(\Sigma, Q, \mu, \sigma) = \inf
_{\hat{f} \in\mathcal
{A}_n}\sup_{f\in\Sigma} E_{f,Q} \llVert
\hat{f} - f \rrVert^2_Q,
\]
where $\mathcal{A}_n$ is the space of all measurable functions of data
to $L_2(Q)$ and $E_{f,Q}$ denotes expectation under the model: $X^i
\sim Q$, $Y^i \mid X^i \sim N(\mu+ f(X^i), \sigma^2)$, independently
across $i = 1, \ldots, n$. When no risk of ambiguity is present, we
shorten the notation $r_n(\Sigma, Q, \mu, \sigma)$ to $r_n$ and call
$r_n$ the minimax risk or the minimax \textit{rate}, when viewed as a
function of the sample size $n$.

Let $\bbN$ denote the set of natural numbers and $\bbN_0 = \bbN\cup
\{0\}$. For any \mbox{$d$-}dimen\-sional multiindex $a = (a_1, \ldots, a_d) \in
\bbN_0^d$ define $\llvert a\rrvert = a_1 + \cdots+ a_d$ and let $D^a$
denote the
mixed partial derivative operator $\partial^{\llvert a\rrvert } /
\partial x_1^{a_1}
\cdots\partial x_d^{a_d}$. For any real number~$b$, let $\lfloor b
\rfloor$ denote the largest integer strictly smaller than $b$. Use the
notation $C^{\alpha, d}$ to denote the Banach space of H\"older
$\alpha$-smooth functions on $[0,1]^d$ equipped~with the norm
\[
\llVert
f\rrVert
_{C^{\alpha, d}} = \sum_{\llvert k\rrvert \le\lfloor\alpha\rfloor}
\bigl\llVert D^kf \bigr\rrVert
_\infty+ \max_{x\neq y \in[0,1]^d} \bigl\llvert D^{\lfloor\alpha\rfloor}f(x)
- D^{\lfloor\alpha\rfloor}(y)\bigr\rrvert / \llVert x - y \rrVert ^{\alpha
- \lfloor\alpha
\rfloor}.
\]
Let $C^{\alpha, d}_1$ denote the unit ball of $C^{\alpha, d}$.

For\vspace*{1.5pt} any $b \in\{0,1\}^p$ and $x = (x_1, \ldots, x_p) \in\bbR^p$,
let $x_b = (x_j\dvtx b_j = 1)$ denote the vector of $\llvert
b\rrvert =\sum_{j=1}^pb_j$
predictors picked by $b$ and let $T^b\dvtx C(\bbR^{\llvert b\rrvert
}) \to C(\bbR^p)$
denote the mapping that takes an $f \in C(\bbR^{\llvert b\rrvert })$
to $T^bf\dvtx x
\mapsto f(x_b)$. Let $\mathcal{B}^{p,d}$ denote the set of all $b \in
\{0,1\}^p$ with $\llvert b\rrvert = d$. We formalize the space of centered,
$p$-variate, $\alpha$-smooth functions of sparsity $d$ and bound
$\lambda$ as
\[
\Sigma^p_S(\lambda, \alpha, d):= \biggl\{\bigcup
_{b \in\mathcal
{B}^{p,d}} T^{b} \bigl(\lambda
C^{\alpha, d}_1\bigr) \biggr\} \cap\mathcal{Z}_p,
\]
where $\mathcal{Z}_p = \{f \in C[0,1]^p\dvtx \int f(x) \,dx=0\}$. The
condition that $f$ is centered can be imposed without any loss of
generality due to the presence of the overall mean parameter $\mu$ in
our regression model. The function spaces $\Sigma_S^p(\lambda, \alpha,
d)$ make up~M2. For M3, we consider additive convolutions of multiple\vspace*{1pt}
$\Sigma^p_S$ spaces with an additional restriction on the number of
components a single predictor can appear in. For $k, \bar d \in\bbN$,
$d \in\bbN^k$ define
\begin{eqnarray*}
\mathcal{B}^{p,k,d,\bar d} &=& \Biggl\{\bigl(b^1, \ldots,
b^k\bigr)\dvtx
b^s \in
\mathcal{B}^{p,d_s}, b^s \ne b^t, 1 \le s \ne t
\le k,
\\
&& \hspace*{121pt}  \sum_{t = 1}^k b^t_j
\le\bar d, 1 \le j \le p
\Biggr\}
\end{eqnarray*}
and
\begin{eqnarray*}
\hspace*{-2pt}&& \Sigma_A^{p,k,\bar d}(\lambda, \alpha, d)
\\
\hspace*{-2pt}&&\qquad = \Biggl\{f =
\sum_{s =
1}^k \lambda_s
T^{b^s} f_s\dvtx
f_s
\in C^{\alpha_s, d_s}_1, 1 \le s \le k,
\bigl(b^1, \ldots, b^k\bigr) \in\mathcal{B}^{p,k,d,\bar d}
\Biggr\} \cap\mathcal{Z}_p.
\end{eqnarray*}
%
In studying minimax rates for a fixed $k$, one can set $\bar d$ as
large as $k$. But in the more interesting large $p$ small $n$ scenario
where $k$ increases with $p$, the use of a fixed $\bar d$ is crucial
for interpreting our results.

For a metric space $(\mathcal{S}, \rho)$, the covering $\varepsilon
$-entropy of a subset $S \subset\mathcal{S}$ is the logarithm of the
minimum number of $\rho$-balls of radius $\varepsilon$ and centers in
$\mathcal{S}$ needed to cover $S$, and is denoted $V(\varepsilon, S,
\rho)$. A finite subset $A \subset S$ is called $\varepsilon$-packing in
$S$ if any two elements of $A$ have a $\rho$-distance at least
$\varepsilon$. The logarithm of the maximal cardinality of an
$\varepsilon
$-packing set in $S$ is called the packing $\varepsilon$-entropy of $S$
and is denoted $C(\varepsilon, S, \rho)$.

\section{Minimax risks for large $p$ small $n$ regression}\label{sec3}
\label{sminimax}
Precise calculations of $r_n$ under M2 and M3 and theoretical results
on whether these rates are achieved in practice are known only under
additional simplifying assumptions on the shape of $f$, or, for
inference tasks that are simpler than prediction. We provide a brief
overview of known results before presenting our main theorem on minimax
$L_2$ risk for regression under M3.

\subsection{A brief overview of existing results}\label{sec3.1}

For linear regression where $\Sigma$ is taken as the set of functions
$f(x) = x^T\beta$ with $\beta$ in an $l_q$ ($q \le1$) ball of $\bbR
^p$ and some additional regularity assumptions are made on the design
matrix, \cite{Raskutti2011} shows that
\begin{eqnarray*}
r_n^2\asymp\cases{ d\log(p/d)/n, &\quad for $q=0$,
\cr
(
\log d/n )^{1-q/2}, &\quad for $q\in(0,1]$,}
\end{eqnarray*}
up to some multiplicative constant, where $d$ is the number of
important predictors. As shown in \cite{Comminges2012}, these rates
are the typical minimax risks associated with variable selection
uncertainty. For $q = 0$, the $l_q$ norm precisely encodes the sparsity
condition of $M_2$. See \cite{Verzelen2012,Wainright2009a}
and \cite{Wainright2009b} for additional results and overviews. Many
authors have established near minimax performance guarantees of various
linear regression methods under the $L_2$ prediction loss; see, for
example, \cite{Bickel2009,candes2007,Meinshausen2009} and~\cite{zhang2008}.

As a nonlinear, nonparametric generalization of the linear model, \cite
{Raskutti2012} considers the completely additive special case of M3
where all $k$ components are univariate and have the same smoothness
$\alpha> 0$ and shows
\[
r_n^2\asymp kn^{-{2\alpha}/({2\alpha+ 1})}+\frac{k\log p}{n}.
\]
Clearly, the minimax risk decomposes into two terms, where the first
term is the sum of minimax risks of estimating each component and the
second term is the variable selection uncertainty.

An entirely different generalization of the linear model is the sparse,
fully nonparametric regression model M2. To the best of our knowledge,
the only minimax rates result in this context is \cite{Comminges2012},
which analyzes minimax risks of support recovery where the objective is
to identify the important predictor rather than estimation of $f$
itself. It is shown that if $d\log(p/d)/n$ is lower bounded by some
positive constant, then for some constant $c>0$,
\[
\inf_{\hat{J}_n}\sup_{f\in\Sigma} P_f(
\hat{J}_n\neq J_f)\geq c,
\]
where $\hat{J}_n$ ranges over all variable selection estimators, that
is, measurable maps of data to the space of all subsets of $\{1, \ldots
, p\}$, $\Sigma$ is the space of all differentiable functions that
depend on only $d$ many predictors and have squared integrable
gradients, and $J_f \subset\{1, \ldots,p\}$ is the index set of truly
important predictors associated with $f$. Because of this result, we
refer to the term $d\log(p/d)/n$ as the risk associated with variable
selection uncertainty. For large $p$, $d\log(p/d)$ is asymptotically
of the same order as the logarithm of ${p \choose d}$, the number of
ways to select $d$ out of $p$ predictors. Any estimation problem
involving high-dimensional variable selection is likely to include a
variable selection uncertainty term $d\log(p/d)/n$ in its minimax rate.

\subsection{New results on minimax rates under \textup{M2} and \textup{M3}}\label{sec3.2}\label{sem2m3}
We calculate minimax $L_2$ risks under the following condition on the
stochastic design:

\renewcommand{\theass}{Q}
\begin{ass}\label{assQ}
$Q$ admits a probability density function (p.d.f.)
$q$ on $[0,1]^p$ such that $\bar q:= \sup_x q(x) < \infty$ and $\inf_{x
\in[1/2-\Delta, 1/2+\Delta]^p} q(x) \ge\underline q$ for some
$\underline q > 0$ and $0 < \Delta\le1/2$.
\end{ass}

The requirement of $q$ being lower bonded on some sub-hypercube inside
$[0,1]^p$ is crucial to obtaining sharp lower bounds on the minimax
risk. This requirement is essentially equivalent to asking that $X$
cannot be reduced to a lower dimension without some loss of
information, for example, $X$ cannot lie on a lower-dimensional
subspace of manifold as assumed under \textup{M1}.


%
%
\begin{theorem}\label{thmM3}
Under Assumption \ref{assQ}, there exist $N_0 \in\bbN$, $0 < \underline C \le
1 \le\overline C$, all depending only on $\bar d$, $\max_s d_s$,
$\min_s \alpha_s$, $\max_s \alpha_s$, $\min_s\lambda_s$, $\max_s
\lambda_s$, such that for all $n > N_0$,
\[
\underline C \underline\varepsilon_n^2 \le
r_n^2 \bigl(\Sigma_A^{p,k,
\bar d}(\lambda,
\alpha, d),Q,\mu,\sigma\bigr) \le\overline C \bar\varepsilon_n^2,
\]
where
\begin{eqnarray*}
\underline\varepsilon_n^2 & =& \sum
_{s=1}^k \lambda_s^2
\biggl(\frac
{\sqrt n \lambda_s } {\sigma} \biggr)^{-4\alpha_s/(2\alpha _s+d_s)}+ \frac{\sigma^2\sum_s d_s}{n} \log
\frac{p}{\sum_s
d_s}
\end{eqnarray*}
and
\begin{eqnarray*}
\bar\varepsilon_n^2 & =& \sum
_{s=1}^k \lambda_s^2
\biggl(\frac
{\sqrt n \lambda_s } {\sigma} \biggr)^{-4\alpha_s/(2\alpha _s+d_s)}+ \frac{\sigma^2\sum_s d_s}{n} \log
\frac{p}{\min_s d_s}.
\end{eqnarray*}
\end{theorem}

%
%
\begin{remark} By choosing $k=1$ and $\bar{d}=1$ in Theorem \ref
{thmM3}, we obtain the minimax risk for M2 as a simple corollary,
%
\begin{equation}
\label{eqminimaxhsn} r_n^2\bigl(\Sigma_S^p(
\lambda, \alpha, d), Q, \mu, \sigma\bigr) \asymp\lambda^2 \biggl(
\frac{\sqrt n \lambda} {\sigma} \biggr)^{-4\alpha/(2\alpha+d)}+\frac{\sigma^2d}{n} \log
\frac{p}d.
\end{equation}
\end{remark}

%
%
\begin{remark}
One can shed light on the scope and limitations of a model by
investigating the conditions needed on the model parameters in order to
bound the model's minimax risk by a given margin. From \eqref
{eqminimaxhsn}, the minimax risk of M2 consists of two terms. The
second term is the typical risk associated with variable selection
uncertainty \cite{Comminges2012} which remains small as long as $\log
p \asymp n^{\beta}$ for some $\beta\in(0,1)$, which gives the
standard large $p$ small $n$ dynamics between sample size and predictor
count. The first term in \eqref{eqminimaxhsn} is the minimax risk of
estimating a $d$-variate, $\alpha$-smooth regression function $f_0$
when there is no variable selection uncertainty. For a fixed smoothness
level $\alpha$, this term remains small as long as $d = o(\log n) =
o(\log\log p)$ under standard large $p$ small $n$ dynamics. In other
words, meaningful statistical learning is possible under M2 only when
the true number of important predictors is much much smaller than the
total predictor count.
\end{remark}

%
%
\begin{remark}
M3 offers a platform to break away from such extreme sparsity
conditions. We consider two special cases for illustration under a
standard large $p$ small $n$ dynamic: $\log p = n^\beta$ for some
$\beta\in(0,1)$, while allowing $k$ to depend on $n$. First, suppose
all additive components $f_s$ have the same dimension ($d_s \equiv d$),
smoothness ($\alpha_s \equiv\alpha$) and magnitude ($\lambda_s
\equiv\lambda$), all of which \mbox{remain} fixed as $k$ increases $n$. This
situation includes as a special case the completely additive framework
of \cite{Raskutti2012}. From Theorem \ref{thmM3}, the associated
minimax risk $r_n^2 \asymp k n^{-{2\alpha}/({2\alpha+ d})} + kd \log
(p / d) / n$ which remains small as long as $k = o(\min\{n^{{2\alpha
}/({2\alpha+ d})}, \log p / n\}) \asymp o(n^\gamma)$ for some $\gamma
\in(0,1)$. Thus, the total number of important predictors, which is of
the order $kd$, could be as large as a fractional power of $\log p$, a
number that is much larger than what is allowed under M2.

In the second case, consider an unbalanced case where $d_s$, $\alpha
_s$ vary with $s$, but remain bounded as $k$ increases with $n$, and
the magnitudes diminish so that the series $\sum_s \lambda_s^{2d_s /
(2\alpha_s + d_s)}$ is convergent. Theorem \ref{thmM3} suggests
that a consistent estimator of $f$ exists in this case as long as $\sum
_{s=1}^kd_s=o(n)$, that is, the total number of important predictors is $o(n)$.
\end{remark}

%
%
\begin{remark}
Consider another unbalanced scenario where $k$ is fixed and one
additive component is much more complex than the rest, that is,
$d_1/\alpha_1\gg d_s/\alpha_s$ for $s=2,\ldots,k$. In this case,
Theorem \ref{thmM3} gives a minimax risk $r_n^2\sim n^{-{2\alpha
_1}/({2\alpha_1+d_1})}+\sum_{s=1}^kd_s\log(p/d_s)/n$, where the
first term is dominated by the largest risk of all additive components,
while the second term is still determined by the overall variable
selection uncertainty. Therefore, the difficulty of estimating a
function with an additive form is determined by the estimation
difficulty of its ``hardest'' component.
\end{remark}

\section{Adaptive near minimax optimality of Bayesian additive
Gaussian process regression}\label{sec4}
\label{sGP}
A Gaussian process (GP) on an Euclidean set $K$ is a random element
$W=(W_x\dvtx x \in\mathcal{X})$ of the supremum-norm Banach space of
continuous functions over $\mathcal{X}$ such that any linear
functional of $W$ is univariate Gaussian \cite{Van2008a}. The
probability law of a GP $W$ is completely determined by the mean and
covariance functions $m(x)=EW_x$ and $\mathcal
{C}(x,x')=E(W_x-m(x))(W_{x'}-m(x'))$ and is denoted by $\gp(m,
\mathcal{C})$. For any function $m\dvtx\mathcal{X}\to\bbR$ and any
nonnegative definite function $\mathcal{C}\dvtx\mathcal{X}\times
\mathcal
{X} \to(0,\infty)$, there exist a GP $W$ with law $\gp(m, \mathcal{C})$.

Adaptivity and near minimax optimality of Bayesian Gaussian process
regression methods are known for low-dimensional applications \cite
{Van2009}. In GP regression, $f$ is assigned a $\gp(m, \mathcal{C})$
prior and inference on $f$ is carried out by summarizing the resulting
posterior distribution given data, which also remains a GP law~\cite
{Rasmussen2006}. Theoretical treatments of GP regression have typically
focused on $m \equiv0$ and $\mathcal{C}(x,x') = \mathcal{C}^{\mathrm
{SE}}(x, x') = \exp(-\llVert x-x' \rrVert ^2)$, the square
exponential covariance
function, with additional hyper-parameters inserted inside the
covariance function \cite{Van2009,Tokdar2010,Choi2007}. In
particular, in order to achieve adaptation to unknown smoothness, \cite
{Van2009} considers as prior distribution the law of a rescaled process
$W^A$ defined as $W^A_x = W_{Ax}$ where $W \sim\gp(0, \mathcal
{C}^{\mathrm{SE}})$ and $A^p$ follows a gamma distribution, and proves the
resulting posterior distribution contracts to the true $f$ at the
minimax rate $n^{-\alpha/(2\alpha+ p)}$ up to a $\log n$ factor when
$f$ is H\"older $\alpha$-smooth. Extensions to anisotropic function
spaces are carried out by \cite{BPD2014}.

\subsection{Additive Gaussian process regression}\label{sec4.1}
\label{ssagp}
For a stochastic process $W = (W_x\dvtx x \in\bbR^p)$, a scalar $a > 0$
and a binary inclusion vector $b \in\{0,1\}^p$, define a \textit
{selective-rescaled} process $W^{a,b} = (W^{a,b}_x\dvtx x \in[0,1]^p)$ by
$W^{a,b}_x = W_{ab\odot x}$ where $\odot$ is the elementwise product
operator. Toward a Bayesian estimation of regression functions $f$
described by M3, we consider the following additive Gaussian process
(add-GP) prior distribution on $f$:
%
\begin{equation}
\label{eqAGP} %
\begin{aligned} f & =& L_1W_1^{A_1,B_1}
+\cdots+ L_KW_K^{A_K, B_K};\qquad K \sim\pi,
\end{aligned}
\end{equation}
where $\pi$ is a probability distribution on $\bbN$, and $L_s
W_s^{A_s, B_s}$ are IID copies of the process $LW^{A,B}$ defined as:
$W\in C(\bbR^p)$, $L \in\bbR_+$ and $(A,B) \in\bbR_+ \times\{0,1\}
^p$ are mutually independent random elements distributed as
%
\begin{eqnarray}\label{eqGPVSp}
W &\sim&\gp\bigl(0,
\mathcal{C}^{\mathrm{SE}}\bigr),\qquad L \sim h,
\nonumber\\[-8pt]\\[-8pt]\nonumber
B &\sim&\Biggl[\bigotimes_{j = 1}^p \operatorname{Be}\biggl(\frac
{1}p\biggr)\Biggr]\Bigg|_{\llvert B\rrvert \le
D_0},\qquad A^{\llvert B\rrvert }
\bigr\rrvert B \sim\operatorname{Ga}(a_1,a_2),
\end{eqnarray}
where $h$ is a density function on $(0,\infty)$ and $a_1, a_2, D_0$
are prespecified, positive valued hyper-parameters.

To complete the add-GP prior specification, we need to specify a prior
distribution on $(\mu, \sigma)$. We consider $(\mu, \sigma) \sim
\pi_\mu\times\pi_\sigma$ where $\pi_\mu$ is a Gaussian
distribution and $\pi_\sigma$ admits density function on $\bbR_+$
with a compact support inside $(0, \infty)$.

\subsection{Posterior contraction rates}\label{sec4.2}

For any $x^{1\dvtx\infty} = (x^1, x^2, \ldots) \in([0,1]^p)^{\infty}$
and any $\theta= (\mu, f, \sigma)$, let $P_\theta(\cdot\mid
x^{1\dvtx\infty})$ denote the conditional distribution of $(Y^i\dvtx
i \in
\bbN)$ given $X^i = x^i$, $i \in\bbN$, under \eqref{eqreg}. Let
$\Pi_n(\cdot\mid(x^i, y^i), 1 \le i \le n)$ denote the posterior
distribution of $\theta$ under the add-GP prior given $(X^i, Y^i) =
(x^i, y^i)$, $1 \le i \le n$. Following \cite{Ghosal2007,Vaart08b},
the posterior contraction rate of the add-GP prior at any $\theta^* =
(\mu^*, f^*, \sigma^{*})$ is said to be at least $\varepsilon_n$ if for
every $x^{1\dvtx\infty}$, other than in a $Q^\infty$-null set,
%
\[
\Pi_n \bigl\{\bigl\llVert\mu+ f - \mu- f^*\bigr\rrVert
_n + \bigl\llvert\sigma- \sigma^*\bigr\rrvert\geq M
\varepsilon_n \mid\bigl(x^i, Y^i\bigr), 1\le i
\le n \bigr\}\stackrel{P_{\theta^*}(\cdot\mid x^{1\dvtx\infty
})} {
\longrightarrow} 0
\]
as $n \to\infty$ for some constant $M$, where ${\llVert \cdot\rrVert}
_n$ denotes
an empirical version of the $L_2(Q)$ norm: $\llVert f\rrVert ^2_n =
(1/n)\sum_{i =
1}^n f^2(x^i)$. It is possible to replace ${\llVert \cdot\rrVert} _n$
with ${\llVert \cdot\rrVert} _Q$ by\vspace*{1pt} appealing to the techniques
developed for GP priors in Section~2.4 in~\cite{Yang2013}, but we omit
the details.

%
%
\begin{theorem}\label{thm2}
Under Assumption \ref{assQ}, for any $\mu^* \in\bbR$, $\sigma^* \in\operatorname
{support}(\pi_\sigma)$ and $f^* \in\Sigma^{p,k,\bar d}_A(\lambda
^*, \alpha^*, d^*)$ with $\max_s d^*_s \le D_0$ and $k \le K_0$, the
posterior contraction rate at $\theta^* = (\mu^*, f^*, \sigma^*)$ is
of the order $\varepsilon_n (\log n)^{(1 + D_0)/2}$ where
\[
\varepsilon_n^2 = \sum_{s=1}^{k}
\lambda^{*2}_s \biggl(\frac{\sqrt
{n}\lambda^*_s}{\sigma^*}
\biggr)^{-{4\alpha^*_{s}}/(2\alpha^*_{s}+d^*_s)}(\log n)^{2q_s}+\frac{\sigma^{*2} \sum_s d^*_s}{n} \log p
\]
with $q_s = (1 + d^*_s) / (2 + d^*_s / \alpha^*_s)$, $1 \le s \le k$,
provided $K_0 \log p \le n\varepsilon_n^2$.
\end{theorem}

When $p$ grows with $n$, add-GP regression essentially employs a
sequence of priors changing with $n$. In this case, it is possible and
useful to also let $K_0$ grow with~$n$ and study posterior contraction
rate at a sequence of $f^* = f^*_n$ changing with~$n$. Theorem \ref
{thm2} remains valid as long as $K_0 \log p \le n\varepsilon_n^2$, the
true number of components $k \le K_0$, $\alpha_s$ are bounded from
above and below and $\max_s \lambda_s$ is bounded.

%
%
\begin{remark}
Related work on estimation of $f$ under M3 includes \cite
{Raskutti2012}, where convergence rates are investigated for an
$M$-estimator with a sparsity penalty on the number of additive
components and smoothness penalties on each components. However, \cite
{Raskutti2012} considers only univariate components. 
In \cite{Taiji2012}, PAC-Bayesian bounds are derived for general
additive regression with additive GP priors. However, \cite{Taiji2012}
assumes that the covariate vector $X$ is pre-divided into $M$ subsets
$(X_{(1)},\ldots,X_{(M)})$ and $f(x)=\sum_{m=1}^Mf_m(x_{(m)})$, with
sparsity constraints on the component functions. Both these studies
assume that important predictors are not shared across components,
which makes the studied methods somewhat restricted in application. A
lack of overlap comes with the technical advantage that $\llVert \sum
_s f_s\rrVert
^2_Q$ decomposes to $\sum_s \llVert f_s\rrVert ^2_Q$ if every $f_s$ has
$Q$-integral 0. In the more general case where components are allowed
to share predictors, a na\"ive application of the Cauchy--Schwarz
inequality gives $\llVert \sum_s f_s\rrVert ^2_Q \le k\sum_s \llVert
f_s^2\rrVert _Q$, but
the multiplication by $k$ results in sub-optimal rates unless $K_0$
grows extremely slowly in $n$.\vadjust{\goodbreak} Our assumption that any predictor can
appear in at most $\bar d$ many components, for some fixed $\bar d$,
overcomes this difficulty with the help of Lemma \ref{lemoverlap}.
\end{remark}

\section{Proofs of the main results}\label{sec5}\label{sear}
\label{sproofs}
\subsection{Minimax rates}\label{sec5.1}

%
%
\begin{lemma}
\label{lemsoln}
For every $\alpha, \lambda$, $d \in\bbN$ there exist $N_0 > 0$, $0
< \underline C \le1 \le\overline C$, such that for any $n > N_0$ and
all $p \in\bbN$, the $\varepsilon_n$ that solves $C(\varepsilon_n,
\Sigma
^p_S(\lambda, \alpha, d), \break {\llVert \cdot\rrVert} ) = n \varepsilon
_n^2/\sigma^2$ satisfies\vspace*{-1pt}
\[
\underline C \le\frac{\varepsilon_n^2}{\lambda^2 ({\sqrt n \lambda
}/{\sigma})^{-{4\alpha}/{(2\alpha+ d)}} + (\sigma^2 / n) \log{p \choose
d}} \le\overline C.
\]
\end{lemma}

\begin{pf}
Let $\varepsilon_1, M_0, M_1$ be as in Lemma \ref{lemsparse}. Without
loss of generality, $M_0 \le1 \le M_1$. Let $\delta_n^2 = \lambda^2
({\sqrt n \lambda}/{\sigma})^{-{4\alpha}/{(2\alpha+ d)}} + (\sigma
^2/{n}) \log{p \choose d}$ and set $N_0$ large enough such that
$\delta_n < \varepsilon_1$ for all $n > N_0$. For the remainder of this
proof, abbreviate $\Sigma^p_S(\lambda, \alpha, d)$ to $\Sigma_S$.
The arguments below mostly rest on the fact that $\varepsilon$-packing
entropy is nonincreasing in $\varepsilon$. Note that
\[
C\bigl(M_1^{1/2}\delta_n, \Sigma_S,
\llVert\cdot\rrVert\bigr) \le C \biggl(\lambda\biggl(\frac{\sqrt n
\lambda}{\sigma}
\biggr)^{-{2\alpha}/(2\alpha+ d)}, \Sigma_S, \llVert\cdot\rrVert
\biggr) \le
M_1 n \delta_n^2 / \sigma^2,
\]
where the second inequality follows by sticking in $\lambda(\sqrt n
\lambda/ \sigma)^{-2\alpha/ (2\alpha+ d)}$ as $\varepsilon$ in Lemma
\ref{lemsparse}. Hence, $\varepsilon_n \le M_1^{1/2} \delta_n$. Also,
by Lemma \ref{lemsparse},
\begin{eqnarray*}
&& C\biggl(\biggl(\max\biggl\{\lambda\biggl(\frac{\sqrt n \lambda}{ \sigma
}\biggr)^{-2\alpha/(2\alpha+ d)},
\frac{\sigma}{\sqrt n} \log^{1/2} \pmatrix{p \cr d} \biggr\}\biggr),
\Sigma_S, \llVert\cdot\rrVert\biggr)
\\[-2pt]
&&\qquad\ge M_0 n \max\biggl\{\lambda\biggl(\frac{\sqrt n \lambda}{
\sigma}
\biggr)^{-2\alpha/(2\alpha+ d)}, \frac{\sigma}{\sqrt n} \log
^{1/2} \pmatrix{p \cr d} \biggr
\}^2 \Big/ \sigma^2
\end{eqnarray*}
and hence $\varepsilon_n^2 \ge M_0 \max\{\lambda(\frac{\sqrt n
\lambda}{ \sigma})^{-2\alpha/(2\alpha+ d)}, \frac{\sigma
}{\sqrt n} \log^{1/2} {p \choose d} \} \ge M_0 \delta_n^2 / 2$.
This proves the result with $\underline C = M_0 / 2$ and $\overline C = M_1$.
\end{pf}

\begin{pf*}{Proof of Theorem \ref{thmM3}}
By Theorem 6 of \cite{Yang1999}, the minimax risk $r_n$ is the
solution to $C(r_n, \Sigma_A^{p,k,\bar d}(\lambda, \alpha, d), {\llVert \cdot\rrVert} _Q) = nr_n^2/\sigma^2$. 
For $1 \le s \le k$, let $\delta_{ns}$ be the solution to $C(\varepsilon
, \Sigma^{p_s}_S(\lambda_s, \alpha_s, d_s), {\llVert \cdot\rrVert} ) =
n\varepsilon
^2/\sigma^2$. From Lemma \ref{lemsoln}, there are $N_s > 0$, $0 <
\underline C_s \le1$, such that\vspace*{2pt} for all $n > N_s$, $\delta_{ns}^2 \ge
\underline C_s\{ \lambda_s^2 ({\sqrt n \lambda_s}/\break {\sigma
})^{-{4\alpha_s}/{(2\alpha_s + d_s)}}\hspace*{-1pt} + (\sigma^2 / n) \log{p_s
\choose d_s}\}$. Denote $\delta_n = (\delta_{n1}, \ldots, \delta
_{nk})$, $n > N = \break \max_s N_s$. Then, by Theorem \ref{thmadd}, with
$b_0 = \underline q^{1/2} \Delta^{\max_s \alpha_s + \max_s d_s/2}$,
%
\begin{eqnarray}\label{eqlowb}
C\biggl(\frac{b_0 \llVert \delta_n \rrVert }2, \Sigma
_A^{p,k,\bar d}(
\lambda, \alpha, d), \llVert\cdot\rrVert_Q\biggr) &\ge&
\frac{1}4\biggl\{\frac{3}4n \frac{\llVert \delta_n \rrVert
^2}{\sigma^2} - k \log2\biggr\}
\nonumber\\[-9pt]\\[-9pt]\nonumber
&\ge&\frac{1}{16} \frac{n\llVert \delta_n \rrVert
^2}{\sigma^2}
\end{eqnarray}
provided $k \log2 \le n\llVert \delta_n \rrVert ^2 / (2\sigma^2)$,
and hence,
\begin{eqnarray*}
r_n^2 & \ge&\frac{\llVert \delta_n \rrVert ^2}{16} \ge\frac{1}{16}
\Biggl\{ \sum_{s = 1}^k \underline
C_s \lambda_s^2 \biggl(\frac{\sqrt n \lambda
_s}{\sigma}
\biggr)^{-4\alpha_s/(2\alpha_s + d_s)} + \frac
{\sigma^2} n \log\pmatrix{p_s \cr d_s} \Biggr\}
\\
& \ge&\underline C \Biggl\{ \sum_{s = 1}^k
\lambda_s^2 \biggl(\frac{\sqrt
n \lambda_s}{\sigma}
\biggr)^{-4\alpha_s/(2\alpha_s + d_s)} + \frac{\sigma^2}n \sum_s
d_s \log p \Biggr\},
\end{eqnarray*}
for some $\underline C$.

Next,\vspace*{1pt} let $\varepsilon_n = (\varepsilon_{n1}, \ldots, \varepsilon_{nk})$
where $\varepsilon_{ns}$ is the solution to $C(\varepsilon_s, \Sigma
_S^p(\lambda_s, \alpha_s, d_s),\break  {{\llVert \cdot\rrVert}} ) = n\varepsilon
_s^2/\sigma
^2$, $1 \le s \le k$. By Lemma \ref{lemsoln}, there are $N_s > 0$,
$\overline C_s \ge1$, such that for all $n > N_s$, $\varepsilon_{ns}^2
\le\overline C_s\{ \lambda_s^2 ({\sqrt n \lambda_s}/{\sigma
})^{-{4\alpha_s}/{(2\alpha_s + d_s)}} + (\sigma^2 / n) \log{p
\choose d_s}\}$. Set $N = \max_s N_s$. By Theorem \ref{thmadd}
again, $C(4 \bar q^{1/2}\sqrt B \llVert \varepsilon_n\rrVert,
\Sigma
_A^{p,k,\bar d}(\lambda, \alpha, d), {\llVert \cdot\rrVert} _Q) \le
n\llVert \varepsilon
_n\rrVert ^2/\sigma^2$, and hence
\[
r_n^2 \le16\bar qB \llVert\varepsilon_n
\rrVert^2 \le\overline C \Biggl\{ \sum_{s = 1}^k
\lambda_s^2 \biggl(\frac{\sqrt n \lambda_s}{\sigma
}
\biggr)^{-4\alpha_s/(2\alpha_s + d_s)} + \frac{\sigma^2}n \sum_s
d_s \log p \Biggr\},
\]
completing the proof.
\end{pf*}

\subsection{Posterior contraction rates of add-GP}\label{sec5.2}
According to \cite{Ghosal2007}, Theorem 1 and Section~7.7, and \cite
{Van2008a}, the conclusion of Theorem \ref{thm2} holds if for
$Q^\infty$-almost every $x^{1\dvtx\infty}$ there exist $\mathcal{F}_n
\subset C(\bbR^p)$, $n \in\bbN$, such that
%
%
\begin{eqnarray}
\Pi\bigl(\bigl\llVert\mu+f-\mu^*-f^*\bigr\rrVert_{n}\leq
\varepsilon_n\bigr)&\geq& e^{-n \varepsilon
_n^2},\label{eq1}
\\
\Pi(\mu+f\notin\mathcal{F}_n)&\leq& e^{-4n \varepsilon_n^2},\label
{eq2}
\\
\log N\bigl(\bar\varepsilon_n,\mathcal{F}_n,\llVert\cdot
\rrVert_{\infty}\bigr)&\leq& n\bar\varepsilon_n^2,\label{eq3}
\end{eqnarray}
where $\bar\varepsilon_n = \varepsilon(\log n)^{(1 + D_0)/2}$ and $\Pi$
denotes the add-GP prior on $(\mu, f, \sigma)$. These conditions map
to one to one to concentration properties of the selective-rescaled
Gaussian processes underlying the add-GP formulation. Without loss of
generality, we assume the prior density $h$ on $L$ is a folded Gaussian
p.d.f., and that \mbox{$a_1 = a_2 = 1$}.

Two important objects associated with any Gaussian process are its
reproducing kernel Hilbert space (RKHS) and concentration function. The
RKHS of any GP $W = (W_x\dvtx x \in\mathcal{X})$, with $\mathcal{X}
\subset\bbR^d$, is defined to be the set $\bbH$ of all function $h\dvtx
\mathcal{X} \to\bbR$ that can be written as $h(x) = EW_x S$ for some
$S$ in the closure of the linear span of the collection of random
variables $\{W_x\dvtx t\in\mathcal{X}\}$ in $L_2$ norm. The set $\bbH$
is a Hilbert space with $\langle EWS_1, EWS_2\rangle_{\bbH} =
ES_1S_2$. With $W$ seen as an element in $C(\mathcal{X})$, its
concentration function at any $w \in C(\mathcal{X})$ is defined as
\[
\phi_w(\varepsilon) = \mathop{\inf}_{h \in\bbH\dvtx \llVert h -
w\rrVert _\infty\le
\varepsilon} \llVert h\rrVert
_{\bbH}^2 - \log\Pi\bigl(\llVert W\rrVert_\infty
\le\varepsilon\bigr),\qquad\varepsilon> 0.
\]
We make use of the following well-known inequalities involving the RKHS
and the concentration function:
%
\begin{eqnarray}
\label{eqb1}  e^{-\phi_w(2\varepsilon)} &\ge&\Pi\bigl(\llVert W -
w\rrVert
_\infty\le2\varepsilon\bigr) \ge e^{-\phi_w(\varepsilon)},
\\
\label{eqb2}  \Pi(W \notin M\bbH_1 + \varepsilon\bbB_1)
&\le& 1 - \Phi\bigl(\Phi^{-1}\bigl(e^{-\phi_0(\varepsilon)}\bigr) + M\bigr),
\\
\label{eqb3}  V\bigl(\varepsilon, M\bbH_1, \llVert\cdot\rrVert
_\infty\bigr) &\le& 1/2 + \phi_0\biggl(\frac{\varepsilon}{2M}
\biggr).
\end{eqnarray}
See Lemma 5.3 of \cite{Vaart08b} for a proof of \eqref{eqb1}. The
inequality \eqref{eqb2} is the well-known Borell's inequality \cite
{Borell75}, and the right-hand side can be further bounded by $\exp\{
-M^2/8\}$ when $M^2/8 \ge\phi_0(\varepsilon)$ since $\Phi^{-1}(u) \ge
-\sqrt{2 \log(1/u)}$ for all $u \in(0,1)$. Inequality \eqref{eqb3}
holds because the right-hand side gives an upper bound to $C(\varepsilon
/(2M), \bbH_1, {\llVert \cdot\rrVert} _\infty)$, since, if $h_1, \ldots
, h_N \in
\bbH_1$ are $\varepsilon/(2M)$-separated in ${\llVert \cdot\rrVert}
_\infty$ then $1
\ge\sum_{j = 1}^N \Pi(W \in h_j + \{\varepsilon/ (2M)\}\bbB_1) \ge N
\exp\{-1/2 - \phi_0(\varepsilon/ (2M))\}$ by \eqref{eqb1}.

For any $b \in\{0,1\}^p$, $a > 0$, let $\bbH^{a,b}$ and $\phi
^{a,b}_w$ denote the RKHS and the concentration function of the
selective-rescaled GP $W^{a,b}$ introduced in Section~\ref{ssagp}. By\vspace*{1pt}
definition, $W^{a,b}$ is isomorphic to a $d$ dimensional, rescaled GP
$\tilde W^a$ with $\tilde W \sim \operatorname{GP}(0, \mathcal{C}^{\mathrm{SE}})$ on
$\bbR
^d$, whose RKHS and concentration function have been studied
extensively in \cite{Van2009}. The following results, which are direct
consequences of Lemmas 4.3, 4.6, 4.7 and 4.8 of \cite{Van2009}, are of
particular interest to us:
%
\begin{eqnarray}
&\displaystyle w \in T^b C^{\alpha,\llvert b\rrvert }\quad\implies\quad
\phi^{a,b}_w(\varepsilon) \le G_0
a^{\llvert b\rrvert } \biggl(\log\frac{a}\varepsilon\biggr)^{1 +
\llvert b\rrvert },&
\nonumber\\[-8pt]\label{eqconc0}  \\[-8pt]
\eqntext{\forall a \ge a_0,
\forall\varepsilon< \varepsilon_0 \wedge{G_1}
{a^{-\alpha}},}
\\
\label{eqcontain}  &\displaystyle a_1^{\llvert b\rrvert /2} \bbH^{a_1,b}_1
\subset a_2^{\llvert b\rrvert /2} \bbH^{a_2,b}_1\qquad \forall0
< a_1 < a_2,&
\\
\label{eqconst}  &\displaystyle h \in\bbH^{a,b}_1\quad\implies\quad\bigl\llvert
h(0)\bigr\rrvert\le1,\qquad \bigl\llVert h - h(0)\bigr\rrVert_\infty\le a
\llvert b\rrvert.&
\end{eqnarray}
In \eqref{eqconc0}, the constants $\varepsilon_0, a_0$, $G_0, G_1$
depend only on $w$ and $\llvert b\rrvert $.

%
%
\begin{lemma}\label{lecp}
Suppose $(\varepsilon_n, n \ge1)$ satisfies $n^{-\gamma_1} \le
\varepsilon
_n \le n^{-\gamma_2}$ for some $0 < \gamma_1 < \gamma_2 < 1/2$ and
$K_0 \log p \le n\varepsilon_n^2$. Then there exists a sequence of sets
$\mathcal{F}_n \subset C[0,1]^d$ satisfying
%
\begin{eqnarray}
\label{eqsieve1} \Pi(\mu+ f \notin\mathcal{F}_n) & \le&\exp\bigl(-4n
\varepsilon_n^2\bigr)
\end{eqnarray}
and
\begin{eqnarray}
\label{eqsieve2}\log N\bigl(\bar\varepsilon_n, \mathcal{F}_n,
\llVert\cdot\rrVert_\infty\bigr) & \le& n\bar\varepsilon_n^2
\end{eqnarray}
with $\bar\varepsilon_n \asymp\varepsilon_n (\log n)^{(1 + D_0)/2}$.
\end{lemma}

\begin{pf}
Let $R_n = K_3 n\varepsilon_n^2$, where $K_3$ is a large constant to be
determined later, and define
\[
\bar L^2_n = R_n,\qquad M_n^2
= 8 K_4 R_n(\log n)^{1 + D_0},\qquad
\delta_n = \frac{\varepsilon_n}{K_0 \llvert b\rrvert \bar L_n M_n},
\]
for some constant $K_4$. By \eqref{eqconc0}, and the fact that $\phi
^{a,b}_0(\varepsilon)$ is nondecreasing in $a$, the constant $K_4$ can be
chosen large enough so that
%
\begin{equation}
\label{eqM2} M^2_n \ge8 \phi_0^{a, b}
\biggl(\frac{\varepsilon_n} {K_0\bar
L_n}\biggr)\qquad\forall b \in\bigcup
_{d \le D_0} \mathcal{B}^{p,d}, \forall a \le
R_n^{1/\llvert b\rrvert },
\end{equation}
for all large $n$. Set $N = \lceil D_0 \log\{R_n / \delta_n\} / (\log
4) \rceil$ and take $\Delta_n(b):= \{\delta_n 4^{j/\llvert b\rrvert
}\dvtx\break  1\le j \le
N\}$. For every $r \in\{0\} \cup\Delta_n(b)$, define
\[
\mathcal{F}^{r,b}_n = \bigcup_{a \in(0, \delta_n] \cup\{r\}
\setminus\{0\}}
2 \bar L_n M_n \bbH^{a,b}_1 +
\frac{\varepsilon
_n}{K_0} \bbB_1.
\]
Consider the sieves
\[
\mathcal{F}_n:= [-\sqrt n, \sqrt n] \oplus\bigcup
_{1 \le k \le K_0} \mathop{\mathop{\bigcup_{b^1,\ldots,b^k \in\bigcup
_{d \le D_0} \mathcal{B}^{p,d},}}_{r_s \in\{0\}\cup\Delta_n(b^s), 1\le
s \le k,}}_
{\sum_s r_s^{\llvert b^s\rrvert } \le4R_n}
\mathcal{F}^{r_1,b^1}_n \oplus\cdots\oplus\mathcal{F}^{r_k,b^k}_n,
\]
for $n \in\bbN$.

Fix any $k \in\{1, \ldots, K_0\}$, $b^1, \ldots, b^k \in\bigcup_{d
\le D_0} \mathcal{B}^{p,d}$ and $a \in\bbR_+^k$ satisfying $\sum_{s
= 1}^k a_s^{\llvert b_s\rrvert } \le R_n$. For $1 \le s \le k$, if
$a_s \le\delta
_n$ set $r_s = 0$, otherwise find $r_s \in\Delta_n(b^s)$ such that
$r_s 4^{-1/\llvert b^s\rrvert } < a_s \le r_s$. Then $\sum_{s = 1}^k
r_s^{\llvert b^s\rrvert } \le
4R_n$ and by \eqref{eqcontain}, $\bar L_n M_n \bbH^{a_s, b^s}_1 +
(\varepsilon_n/K_0) \bbB_1 \subset\mathcal{F}^{r_s,b^s}_n$ for all $1
\le s \le k$. Therefore,
\begin{eqnarray*}
&& \Pi\bigl\{\mu+f \notin\mathcal{F}_n \mid K = k,
\bigl(A_s, B^s\bigr) = \bigl(a_s,
b^s\bigr), 1 \le s \le k\bigr\}
\\
&&\qquad  \le\Pi\bigl(\llvert\mu\rrvert> \sqrt n\bigr) + \sum
_{s = 1}^k \Pi\biggl\{L_s
W^{a_s,b^s} \notin\bar L_n M_n\bbH^{a_s,b^s}_1
+ \frac{\varepsilon
_n}{K_0} \bbB_1 \biggr\}
\\
&&\qquad  \le e^{-n/2} + \sum_{s = 1}^k
\biggl[\Pi(L_s > \bar L_n) + \Pi\biggl\{W^{a_s,b^s}
\notin M_n\bbH^{a_s,b^s}_1 + \frac{\varepsilon
_n}{K_0\bar L_n}
\bbB_1 \biggr\} \biggr]
\\
&&\qquad \le e^{-n/2} + k \bigl\{e^{-\bar L^2_n} + e^{-M^2_n / 8}\bigr\}
\\
&&\qquad \le3ke^{-R_n}
\end{eqnarray*}
for all large $n$, by \eqref{eqb2} and \eqref{eqM2} and the fact
$R_n = o(n)$. Consequently,
\begin{eqnarray*}
\hspace*{-4pt}&& \Pi(\mu+f \notin\mathcal{F}_n)
\\
\hspace*{-4pt}&&\qquad \le \mathop{\max
_{1 \le k \le K_0,}}_{b^1, \ldots, b^k \in\bigcup_{d \le D_0} \mathcal
{B}^{p,d}} \lleft\{\Pi\lleft(\sum
_{s = 1}^k A_s^{\llvert b^s\rrvert } >
R_n \bigg|
B^s=b^s,
s = 1, \ldots,k
\rright) + 3ke^{-R_n}
\rright\}
\\
\hspace*{-4pt}&&\qquad  \le \Pi(G > R_n) + 3K_0 e^{-R_n}
\end{eqnarray*}
with $G \sim \operatorname{Ga}(K_0, 1)$. Notice $\Pi(G > R_n) \le\exp\{
-R_n/2 + K_0 \log2\}$. Therefore, by the assumption on $K_0$, $\Pi
(\mu+ f \notin\mathcal{F}_n)$ is bounded by $\exp(-4n\varepsilon
_n^2)$ for all large $n$, provided $K_3$ is chosen suitably large.

By \eqref{eqconst}, when $r = 0$, $\mathcal{F}^{r,b}_n \subset2
\bar L_n M_n \cdot[-1,1] + (2\varepsilon_n / K_0) \bbB_1$, and hence
could be covered by $\lceil4\bar L_n M_nK_0 / \varepsilon_n\rceil$ many
or fewer balls of supremum norm radius $3\varepsilon_n / K_0$. When $r >
0$, by \eqref{eqb3}, at most another $1/2 + \phi_0^{r,b}(\varepsilon_n
/\break  (2\bar L_n M_nK_0))$ many balls\vspace*{1pt} may be needed to maintain
$3\varepsilon
_n / K_0$ covering. Therefore, by \eqref{eqconc0}, $V(3\varepsilon_n /
K_0, \mathcal{F}^{r,b}_n, {\llVert \cdot\rrVert} _\infty) \le D_1 \{
r^{\llvert b\rrvert } (\log
n)^{1 + D_0} + \log n\}$ for every $r \in\{0\} \cup\Delta_n(b)$, for
some constant $D_1$ that depends only on $D_0$, as long as $\llvert
b\rrvert \le
D_0$. Consequently,
\[
V\bigl(4\varepsilon_n, \mathcal{F}_n, \llVert\cdot\rrVert
_\infty\bigr) \le D_1 \bigl\{R_n (\log
n)^{1 + D_0} + K_0 \log n\bigr\} + \log M,
\]
where $M$ is the size of the finite set $\{((r_1, b^1), \ldots, (r_k,
b^k))\dvtx 1 \le k \le K_0, b^s \in\bigcup_{d \le D_0} \mathcal
{B}^{p,d},r_s \in\{0\}\cup\Delta_n(b^s), 1\le s \le k\}$. This
proves the result because of the assumption on $K_0$, since $\log M \le
\log[K_0 \{p^{D_0}(N + 1)\}^{K_0}] \le C_6 K_0 \log p$ for some
constant $C_6$ that depends only on $D_0$.
\end{pf}

%
%
\begin{lemma}\label{lepc}
Under the conditions of Theorem \ref{thm2}, for $Q$-almost every
$x^{1\dvtx\infty}$, $\Pi(\llVert \mu+ f -\mu^* - f^*\rrVert _n \le
\varepsilon_n) \ge
\exp(-n\varepsilon_n^2)$ for all large $n$ where
\[
\varepsilon_n^2 \asymp\sum_{s = 1}^k
\lambda^{2}_s \biggl(\frac{\sqrt
n \lambda_s}{\sigma^*}
\biggr)^{-4\alpha_s/(2\alpha_s +
d_s)} (\log n)^{2q_s} + \frac{\sigma^{*2}\sum_{s} d_s}{n}\log p
\]
with $q_s = (1 + d_s) / (2 + d_s / \alpha_s)$, $s = 1, \ldots, k$.
\end{lemma}

\begin{pf}
By Lemma \ref{lenormc}, with $\mathcal{F}_n$ as in \eqref
{eqsieve2}, we only need to show $\Pi(\llVert \mu+ f -\mu^* -
f^*\rrVert _Q \le
\varepsilon_n, \mu+f\in\mathcal{F}_n, \llVert \mu+f-\mu^*-f^*\rrVert
_{\infty
}\leq1) \ge\exp(-n\varepsilon_n^2)$.\vspace*{1pt}
By inequality \eqref{eqsieve1} and the fact ${\llVert \cdot\rrVert} _Q
\le
\bar q^{1/2}{\llVert \cdot\rrVert} $ it suffices to show that
\[
\Pi\bigl(\bigl\llVert\mu+ f -\mu^* - f^*\bigr\rrVert\le\varepsilon_n,
\bigl\llVert\mu+f-\mu^*-f^*\bigr\rrVert_{\infty}\leq1\bigr) \ge\exp
\bigl(-n\varepsilon_n^2\bigr).
\]
%
We can write $f^* = \sum_{s = 1}^k \lambda_s T^{b^s} f^*_s$ where
$b^s \in\mathcal{B}^{p, d_s}$, $f^*_s \in C^{\alpha_s, d_s}_1 \cap
\mathcal{Z}_d$, $1 \le s \le k$ and $\max_{1 \le j \le p} \sum_{s =
1}^k b^s_j \le\bar d$. Let $\delta_{ns} = \lambda_s (\sqrt n \lambda
_s / \sigma)^{-2\alpha_s / (4\alpha_s + d_s)} (\log n)^{q_s}$, $1
\le s \le k$ and $\delta_n = (\delta_{n1}, \ldots, \delta_{nk})$.
Set $B = 1 + \max_s d_s (\bar d - 1)$.

For any $a > 0$, $b \in\{0,1\}^p$ define the Gaussian variable
$U^{a,b} = \int W^{a,b}_x \,dx$. Then the Gaussian process $V^{a,b} =
W^{a,b} - E(W^{a,b} \mid U^{a,b})$ satisfies $\int V^{a,b}_x \,dx = 0$, and
is independent of the process $ E(W^{a,b} \mid U^{a,b}) = Z \psi^{a,b}$
where $Z \sim N(0,1)$ and $\psi^{a,b}(x) = \operatorname{cov}(U^{a,b}, W^{a,b}_x) /
\operatorname{var}^{1/2}(U^{a,b})$, $x \in[0,1]^p$. By\vspace*{2pt} Cauchy--Schwarz inequality,
$\llVert
\psi^{a.b}\rrVert _\infty\le1$. Clearly,\vspace*{2pt} $W^{a,b}$ decomposes as $W^{a,b}
= V^{a,b} + Z \psi^{a,b}$.

Therefore, for any $\ell, a \in\bbR_+^k$ and given $K = k$, $(L_s,
A_s, B^s) = (\ell_s, a_s, b^s)$, $1 \le s \le k$, we can decompose the
additive-GP process $f$ as $f = \sum_{s = 1}^k \ell_s\times\break  Z_s \psi^{a_s,
b^s} + \sum_{s = 1}^k \bar f_s$, where $Z_s$ are independent $N(0,1)$
variables, $\bar f_s$ are mutually independent with probability laws
same as those of $\ell_s V^{a_s, b_s}$, and these two sets of random
quantities are independent. Consequently, for large enough $n$,
\begin{eqnarray*}
\hspace*{-4pt}&& \Pi\bigl\{\bigl\llVert f - f^*\bigr\rrVert\le\sqrt{1 + 25B}
\llVert\delta_n\|,\bigl\| f-f^*\bigr\|_{\infty}\leq1/2
\\
\hspace*{-4pt}&&\hspace*{26pt} \mid K = k,\bigl(L_s, A_s,
B^s\bigr) = \bigl(\ell_s, a_s,
b^s\bigr), 1 \le s \le k \bigr\}
\\
\hspace*{-4pt}&&\qquad  \ge\Pi\biggl(\biggl\llVert{ \sum_s}
\ell_s Z_s \psi^{a_s, b^s}\biggr\rrVert\le\llVert
\delta_n \rrVert, \biggl\llVert{ \sum_s}
\ell_s Z_s \psi^{a_s, b^s}\biggr\rrVert
_{\infty}\leq1/4 \biggr)
\\
\hspace*{-4pt}&&\quad\qquad{}\times\Pi\biggl\{\biggl\llVert{ \sum_s}
\bigl(\bar f_s - \lambda_s T^{b^s}
f^*_s\bigr)\biggr\rrVert\le5\sqrt B \llVert\delta_n
\rrVert,
\biggl\llVert{ \sum_s} \bigl(\bar
f_s - \lambda_s T^{b^s} f^*_s
\bigr)\biggr\rrVert_{\infty} \leq1/4\bigg|
\\
\hspace*{-4pt}&& \hspace*{157pt}
K =
k,\bigl(L_s, A_s, B^s\bigr) =
\bigl(\ell_s, a_s, b^s\bigr), 1 \le s \le k
\biggr\}
\\
\hspace*{-4pt}&&\qquad \ge\Pi\biggl(\biggl\llVert{ \sum_s}
\ell_s Z_s \psi^{a_s, b^s}\biggr\rrVert\le\llVert
\delta_n \rrVert, \biggl\llVert{ \sum_s}
\ell_s Z_s \psi^{a_s, b^s}\biggr\rrVert
_{\infty} \leq1/4 \biggr)
\\
\hspace*{-4pt}&&\quad\qquad{}\times\prod_{s = 1}^k \Pi\bigl(\bigl
\llVert\ell_s V^{a_s, b^s} - \lambda_s
T^{b^s} f^*_s\bigr\rrVert\le5\delta_{ns},\bigl
\llVert\ell_s V^{a_s, b^s} - \lambda_s
T^{b^s} f^*_s\bigr\rrVert_{\infty} \leq
\delta_{ns}\bigr),
\end{eqnarray*}
because of Lemma \ref{lemoverlap}, since by the assumption on $f^*_s$
and the construction of $\bar f_s$, we have for every $1 \le s \le k$,
$\langle\bar f_s - \lambda_s T^{b^s} f^*_s, \bar f_t - \lambda_s
T^{b^s} f^*_s\rangle_R \ne0$ for at most $r_s = 1 + d_s (\bar d - 1)$
many $1 \le t \le k$.

If $\ell_s \in\lambda_s \cdot[1, 1 + \delta_{ns}]$, then $\{\llVert
\ell_s V^{a_s, b^s} - \lambda_s T^{b^s} f^*_s\rrVert \le5\delta_{ns}\}
\supset\{\lambda_s \llVert V^{a_s, b^s} - T^{b^s} f^*_s\rrVert _\infty
\le
4\delta_{ns}\}\supset\{\lambda_s \llVert W^{a_s, b^s} -
T^{b^s} f^*_s\rrVert
_\infty\le2\delta_{ns}\}$. When $a^{d_s}_s \in(G_1 /\break  \delta
_{ns})^{d_s/\alpha_s} \cdot[1, 2]$, where $G_1$ is as in \eqref
{eqconc0}, the last probability can be lower bounded by $\exp\{
-G_2(\lambda_s / \delta_{ns})^{d_s / \alpha_s} \log(1 / \delta
_{ns})^{1 + d_s}\} \ge\exp\{-G_2 n\delta_{ns}^2 / \sigma^2\}$ for
some constant $G_2$, for all large $n$, by \eqref{eqb1} and \eqref
{eqconc0}. For the same choices of $\ell_s, a_s$, $1 \le s \le k$,
$\Pi(\llVert \sum_s \ell_s Z_s \psi^{a_s, b^s}\rrVert \le\llVert
\delta_n \rrVert,\llVert \sum_s \ell_s Z_s \psi^{a_s, b^s}\rrVert
_{\infty} \le1/4) \ge\exp\{-G_3  k \log n\}$ for some constant $G_3$, for all
large $n$. Therefore, by the
assumption on $K_0$,
\begin{eqnarray*}
&& \Pi\bigl(\bigl\llVert f - f^*\bigr\rrVert\le\sqrt{1 + 25 B}\llVert
\delta_n \rrVert,\bigl\llVert f-f^*\bigr\rrVert_{\infty
}
\le1/2\bigr)
\\
&&\qquad \ge\exp\bigl(-G_4 2n\llVert\delta_n \rrVert
^2/ \sigma^2\bigr)\Pi(K = k)
\\
&&\quad\qquad{}\times \prod
_{s =
1}^k \bigl\{\Pi\bigl(L_s \in
\lambda_s \cdot[1, 1 + \delta_{ns}]\bigr)\Pi\bigl(B^s = b^s\bigr)
\\
&&\hspace*{62pt}{} \times\Pi
\bigl(A^{d_s}_s \in({G_1}/{
\delta_{ns}})^{d_s/\alpha_s} \cdot[1, 2]\mid\llvert B_s
\rrvert= d_s\bigr) \bigr\}
\\
&&\qquad \ge G_5 \exp\biggl\{-G_6 n\biggl\{\llVert
\delta_n \rrVert^2 + \frac{\sigma^2 \sum_s
d_s}{n} \log p\biggr\}
\biggr\}
\end{eqnarray*}
for all large $n$ for some constants $G_5, G_6$ that depend only on
$\max_s d_s$, $\min_s \lambda_s$, $\max_s \lambda_s$, $\min_s
\alpha_s$ and $\max_s \alpha_s$. This proves the result since $\mu$
is independent of $f$ and $\Pi(\llvert \mu- \mu^*\rrvert \le\min\{
\llVert \delta_n \rrVert,1/2\}) \ge\exp\{- G_7 \log n\}$ for some
constant $G_7$.
\end{pf}

\begin{pf*}{Proof of Theorem \ref{thm2}}
Equations~\eqref{eq1}--\eqref{eq3} are implied by Lemmas \ref{lepc}~and~\ref{lecp} with the $\varepsilon_n$ given in Theorem \ref{thm2}.
\end{pf*}

\section{Auxiliary results}\label{sec6}
\label{saux} In this section, we provide a number of auxiliary
results on packing and covering entropies of regular, sparse and
additive H\"older spaces.

%
%
\begin{lemma}
\label{lembasic}
For every $\alpha> 0$, $d \in\bbN$ there exist $\varepsilon_0 > 0$,
$M_0 > 0$ such that for all $\varepsilon< \varepsilon_0$ there are $N
\ge
\exp\{M_0 (1/\varepsilon)^{d/\alpha}\}$ functions $f_0, \ldots, f_N
\in C^\infty(\bbR^d)$ satisfying $f_0 \equiv0$ and
%
\begin{eqnarray}
\label{eqfsupp}
\operatorname{support}(f_i) &\subset&[0,1]^d,\qquad f_i
\mid_{[0,1]^d} \in C^{\alpha,d}_1,\qquad0 \leq i \leq N,
\\
\label{eqfzero}  \qquad \int_\bbR f_i(u_1,
\ldots, u_d) \,du_j &=& 0,\qquad0 \le i \le N,1 \le j \le d,
\\
\label{eqfsep}\llVert f_i - f_k \rrVert&\ge&\varepsilon, \qquad 0
\le i < k \le N.
\end{eqnarray}
\end{lemma}

\begin{pf}
Our proof follows the calculations in \cite{Tsybakov2009},
Section~2.6.2, suitably adapted to handle $L_2$ norm and
condition \eqref{eqfzero}. Let $\mathcal{K} \in C^\infty(\bbR^d)$
such that
%
\begin{equation}
\qquad\operatorname{support}(\mathcal{K} ) = [-1,1]^d,\qquad \int\mathcal{K}
(u_1, \ldots, u_d) \,du_j = 0, j = 1, \ldots,
d.
\end{equation}
For example, one could take $\mathcal{K} (x_1,\ldots,x_d) = \prod_{j
= 1}^d \mathcal{K} _0(x_j)$ where $\mathcal{K} _0(t) =
te^{-{1}/({1-t^2})}I(\llvert t\rrvert \leq1)$, $t \in\bbR$.

Fix an arbitrary $h \in(0,1/2)$ and take $m= \lceil{1}/({2h})
\rceil$, $M=m^d$
and a rectangular grid $\{x^k\dvtx k = 1, \ldots, M\}$ on $[0,1]^d$
consisting of the $M$ grid points $(\frac{j_1-1/2}{m},\ldots,\frac
{j_d-1/2}{m})$, $(j_1,\ldots,j_d)\in\{1,\ldots,m\}^d$. We assume $h$
is small enough so that $M \ge8$. For each $1\le k \le M$, the
function $\phi_k$ defined as
%
\begin{equation}
\label{eqphi} \phi_k(x)= \frac{1}{\llVert \mathcal{K}\rrVert
_{C^{\alpha, d}}}h^{\alpha
}
\mathcal{K} \biggl(\frac{x-x^k}{h} \biggr),\qquad x\in[0,1]^d
\end{equation}
has support inside $x^k + [-h,h]^d$ and belongs to $C^{\alpha, d}_1$.
Let $\Omega= \{0,1\}^M$ and for each $\omega\in\Omega$ define
$f_\omega= \sum_{k = 1}^M \omega_k \phi_k$. Clearly, each $f_\omega
$ is supported on $[0,1]^d$ and $\int f_\omega(u_1, \ldots, u_d)\,du_j
= 0$ for every $j =1,\ldots,d$. Also, since $\phi_k$'s are shifted
copies of each other with disjoint supports, each $f_\omega\in
C^{\alpha, d}_1$ and
%
\begin{eqnarray}
\label{eqdiff} \llVert f_{\omega} - f_{\omega'}\rrVert &=& \Biggl\{
\sum_{k=1}^M\bigl(\omega_k-
\omega'_k\bigr)^2\int\phi_k^2(x)\,dx
\Biggr\}^{1/2}
\nonumber\\[-8pt]\\[-8pt]\nonumber
&=& h^{\alpha+d/2}\frac{\llVert \mathcal{K}
\rrVert }{\llVert \mathcal{K} \rrVert _{C^{\alpha,d}}} \rho
^{1/2}\bigl(\omega, \omega'\bigr),
\end{eqnarray}
where $\rho(\omega, \omega')=\sum_{k=1}^MI(\omega_k\neq\tilde
\omega_k)$ denotes the Hamming distance.

By the Varshamov--Gilbert bound \cite{Tsybakov2009}, Lemma 2.9, there
are $N \ge2^{M/8}$ binary strings $\omega^{(0)},\ldots,\omega^{(N)}
\in\Omega$, with $\omega^{(0)} = 0$, satisfying $\rho(\omega
^{(k)},\omega^{(k')})\geq{M}/{8}$, $0\leq k < k' \leq N$. Then $f_i:=
f_{\omega^{(i)}}$, $0 \le i \le N$, satisfy \eqref
{eqfsupp}--\eqref{eqfzero} and
\[
\llVert f_i - f_k \rrVert\ge h^{\alpha+ d/2}
\frac{\llVert \mathcal{K} \rrVert }{\llVert
\mathcal{K} \rrVert _{C^{\alpha, d}}} \sqrt{\frac{M}8} \ge M_1
h^{\alpha
},\qquad1 \le i < k \le N,
\]
where $M_1 = \llVert \mathcal{K} \rrVert / \{2^{(d + 3)/2} \llVert
\mathcal{K} \rrVert
_{C^{\alpha, d}}\}$ depends on only $\alpha$
and $d$. This proves
the result since with $\varepsilon= M_1h^\alpha$, which could be
arbitrarily small, we get $N \ge\exp\{M (\log2) / 8\} \ge\exp\{M_0
(1/\varepsilon)^{d/\alpha}\}$ where $M_0 = (M_1^{d/\alpha} \log2) /
2^{d + 3}$ depends on only $d$ and $\alpha$.
\end{pf}

%
%
\begin{lemma}
\label{lemsparse}
For every $\alpha, L > 0$, $d \in\bbN$ there exist $\varepsilon_1,
M_0, M_1 > 0$ such that for any $\varepsilon< \varepsilon_1$ and all $p
\in\bbN$
\[
M_0 (L/\varepsilon)^{d/\alpha} + \log\pmatrix{p \cr d} \le C\bigl(
\varepsilon, \Sigma_S^p(\lambda, \alpha, d), \llVert\cdot
\rrVert\bigr) \le M_1 (L/\varepsilon)^{d/\alpha} + \log\pmatrix{p \cr d},
\]
and, an $\varepsilon$-packing set satisfying the above lower bound may be
obtained entirely with $C^\infty(\bbR^p)$ functions.
\end{lemma}

\begin{pf} It suffices to prove for $L = 1$ since $C(\varepsilon, L\Sigma,
{\llVert \cdot\rrVert} ) = C(\varepsilon/L, \Sigma, \llVert \cdot
\rrVert )$ for any set
$\Sigma$. By Lemma \ref{lembasic} there exist $\varepsilon_0, M_0$
such that for any $\varepsilon< \varepsilon_0$ there are functions $f_0
\equiv0, f_1, \ldots, f_N \in C^\infty(\bbR^d)$ satisfying \eqref
{eqfsupp}--\eqref{eqfsep} with $\log N \ge M_0 (1/\varepsilon
)^{d/\alpha}$. Therefore, the set
%
\begin{equation}
\label{eqtset} \mathcal{T}^{\alpha, d,p}(\varepsilon) = \mathop{\bigcup
_{b\in\{0,1\}^p}}_{\llvert b\rrvert = d} \bigl\{T^b f_i\dvtx 1 \le i
\le N\bigr\}
\end{equation}
is a subset of $\Sigma^p_S(1, \alpha, d)$. By \eqref{eqfzero}, for
any $b \ne b' \in\{0,1\}^p$, $\langle T^b f_i, T^{b'} f_k\rangle= 0$
for all $1 \le i, k \le N$. Hence, $\mathcal{T}^{\alpha, d,
p}(\varepsilon)$ is $\varepsilon$-separated in ${\llVert \cdot\rrVert}
$ since $\llVert T^b
f_i - T^{b'} f_k\rrVert = \llVert f_i - f_k \rrVert \ge\varepsilon$
by\vspace*{1pt} \eqref{eqfsep}
if $b = b'$ and $\llVert T^b f_i - T^{b'} f_k\rrVert = \llVert f_i
\rrVert + \llVert f_k\rrVert \ge
\varepsilon$ by \eqref{eqfsep} and the fact that $f_0 \equiv0$. This
gives the lower bound on $C(\varepsilon,\break  \Sigma_S^{p}(1, \alpha, d),
{\llVert \cdot\rrVert} )$ since the cardinality of $\mathcal{T}^{\alpha
, d, p}$ is
${p \choose d} N$.

It is well known that for every $\alpha> 0$, $d \in\bbN$ there exist
$\varepsilon_0', M_0' > 0$ such that for all $\varepsilon< \varepsilon_0'$,
$V(\varepsilon, C^{\alpha, d}_1, {\llVert \cdot\rrVert} ) \le M_0'
(1/\varepsilon
)^{d/\alpha}$ \cite{Tsybakov2009}, Section~2.6.1, and \cite
{Lorentz1996}.
Since a union of sets is covered by the union of their covers, it
follows that $V(\varepsilon, \Sigma_S^p(L\lambda, \alpha, d), {\llVert \cdot\rrVert} ) \le M_0' (1 / \varepsilon)^{d/\alpha} + \log{p \choose
d}$ for all
$0 < \varepsilon< \varepsilon_0'$. Consequent\-ly, $C(\varepsilon, \Sigma
^p_S(\lambda, \alpha, d), {\llVert \cdot\rrVert} ) \le V(\varepsilon/
2, \Sigma
^p_S(\lambda, \alpha, d), {\llVert \cdot\rrVert} ) \le M_0'2^{d/\alpha}
(1/\varepsilon)^{d/\alpha} + \log{p \choose d}$ for all $\varepsilon<
\varepsilon_0'$. This proves the result with $M_1 = M_0'2^{d/\alpha}$
and $\varepsilon_1 = \min(\varepsilon_0, \varepsilon_0')$.
\end{pf}

%
%
\begin{lemma}
\label{lemadd}
Let $\bbH_1, \ldots, \bbH_k$ be mutually orthogonal subsets of a
Hilbert space $(\bbH, {\llVert \cdot\rrVert} _{\bbH})$. Then, for any
$\delta\in
\bbR_+^k$ and $c \in(0,1)$
\[
C\Biggl(c\llVert\delta\rrVert, \bigoplus_{s=1}^k
\bbH_s, \llVert\cdot\rrVert_{\bbH}\Biggr) \ge
\frac{1}4 \Biggl\{\frac{1 - c^{2}}{C^*}\sum_{s=1}^k
C\bigl(\delta_s, \bbH_s, \llVert\cdot\rrVert
_{\bbH}\bigr) - k \log2 \Biggr\},
\]
where $C^*:= \sup_{1 \le s,t \le k} \{\delta_s^{-2}C(\delta_s, \bbH
_s, {\llVert \cdot\rrVert} _{\bbH})\}/\{\delta_t^{-2}C(\delta_t, \bbH
_t, {\llVert \cdot\rrVert} _{\bbH})\}$.
\end{lemma}

\begin{pf}
For every $1 \le s \le k$, let $\mathcal{H}_s$ denote a maximal
$\delta_s$-packing set of $\bbH_s$ with $C_s:= \log\llvert \mathcal
{H}_s\rrvert
= C(\delta_s, \bbH_s, {\llVert \cdot\rrVert} _{\bbH})$. Take $\Omega=
\mathcal
{H}_1 \times\cdots\times\mathcal{H}_k$ and let $F = (F_1, \ldots,
F_k)$ be a random element in $\Omega$ with the uniform probability
distribution. Fix an $M \in\bbN$ such that
\[
\frac{1}2 \biggl\{\frac{1 - c^2}{C^*} \sum_s
C_s - k \log2 \biggr\} < 2\log M < \frac{1 - c^2}{C^*} \sum
_s C_s - k \log2,
\]
and let $F^j$, $j = 1, \ldots, M$, be IID copies of $F$. If
%
\begin{equation}
\label{eqplow} P\biggl\{\biggl\llVert\sum_s
F^i_s - \sum_s
F^j_s\biggr\rrVert\ge c \llVert\delta\rrVert,
\forall1 \le i < j \le M\biggr\} > 0,
\end{equation}
then $\Omega$ contains a subset $\Omega_0$ with at least $M$ elements
such that for any two \mbox{$f, f' \in\Omega$}, $\llVert \sum_s f_s - \sum_s
f'_s\rrVert > c\llVert \delta\rrVert $. This would prove the result.

The probability value in \eqref{eqplow} is at least $1- M(M-1)/2\cdot
P\{\llVert \sum_s F^{1}_s - \sum_s F^{2}_s\rrVert < c\llVert \delta
\rrVert \}$, and hence\vspace*{1pt}
it suffices to show $P\{\llVert \sum_s F^{1}_s - \sum_s
F^{2}_s\rrVert <c\llVert \delta\rrVert \} \le1/M^{2}
$. Define $Z_s = I(F^1_s \ne F^2_s)$, $s = 1, \ldots, k$, which are
independent binary variables with $Z_s \sim\operatorname
{Bernoulli}(1-e^{-C_s})$. By orthogonality of $\bbH_1, \ldots, \bbH_k$,
\[
\biggl\llVert\sum_s F^{1}_s
- \sum_s F^{2}_s\biggr
\rrVert^2=\sum_{s=1}^k\bigl
\llVert F^{1}_s-F^{2}_s\bigr
\rrVert^2\geq\sum_{s=1}^k
\delta_s^2Z_s,
\]
and hence it suffices to show
%
\begin{equation}
\label{eqplow2} P\biggl( \sum_{s}
\delta_s^2 Z_s < c^2\llVert
\delta\rrVert^2\biggr) \le1 / M^2.
\end{equation}
By Markov's inequality, for any $\lambda> 0$,
\begin{eqnarray*}
P \biggl( \sum_s \delta_s^2Z_s
 < c^2\llVert\delta\rrVert^2 \biggr) &\leq& P \bigl
\{e^{-\lambda\sum_s \delta_s^2Z_s}>e^{-\lambda c^2\llVert
\delta\rrVert ^2} \bigr\}
\\
& \leq& e^{\lambda c^2\llVert \delta\rrVert ^2}\prod_{s=1}^kE
\bigl\{e^{-\lambda
\delta_s^2Z_s} \bigr\}
\\
&\leq& e^{\lambda c^2\llVert \delta\rrVert
^2}\prod
_{s=1}^k \bigl\{e^{-C_s}
+e^{-\lambda\delta_s^2} \bigr\}
\\
& =&e^{-\lambda(1 - c^2)\llVert \delta\rrVert ^2}\prod_{s=1}^k \bigl
\{1+e^{\lambda
\delta_s^2-C_s} \bigr\}.
\end{eqnarray*}
By\vspace*{1pt} the assumption on $\delta$, ${C_s \delta_t^2 }/({\delta_s^2 C^*})
\le C_t \le{C^* C_s \delta_t^2 }/({ \delta_s^2})$ for every $1 \le
s, t \le k$, and hence, $\delta_s^2 \le C^* C_s \llVert \delta\rrVert
^2 / \sum_t C_t \le C_s / \lambda$ when we set $\lambda= \sum_s C_s /\break 
(\llVert
\delta\rrVert ^2 C^*)$. Consequently,
\[
P \biggl(\sum_s \delta_s^2
Z_s < c^2\llVert\delta\rrVert^2 \biggr)
\le2^k e^{-\lambda(1 - c^2)\llVert \delta\rrVert ^2} = e^{-(1 - c^2)
\sum_s C_s
/ C^* + k \log2} \le1/M^2,
\]
which completes the proof.
\end{pf}

%
%
\begin{theorem}
\label{thmadd}
Suppose $k \max_s d_s \le p$ and set $p_s = \lfloor pd_s / \sum_t d_t
\rfloor$, $1 \le s \le k$. Under Assumption \ref{assQ}, for any $\delta\in
\bbR_+^k$,
\begin{eqnarray*}
&& C \bigl({\underline{q}{}^{1/2} \Delta^{\max_s (\alpha_s +
d_s/2)}\llVert\delta
\rrVert}/{2},\Sigma_A^{p, k, \bar d}(\lambda, \alpha, d), {\llVert \cdot\rrVert}_Q \bigr)
\\
&&\qquad \ge\frac{1}4 \Biggl\{\frac{3}{4C^*} \sum
_{s = 1}^k C\bigl(\delta_s,
\Sigma^{p_s}_S(\lambda_s,
\alpha_s, d_s), \llVert\cdot\rrVert\bigr) - k \log2
\Biggr\},
\\
&& C\bigl(4\bar{q}{}^{1/2}\sqrt{B}\llVert\delta\rrVert,
\Sigma_A^{p, k, \bar
d}(\lambda, \alpha, d), \llVert\cdot\rrVert
_Q\bigr) \le\sum_{s = 1}^k C
\bigl(\delta_s, \Sigma^p_S(
\lambda_s, \alpha_s, d_s), \llVert\cdot
\rrVert\bigr)
\end{eqnarray*}
with $C^* = \sup_{1 \le s,t \le k} \{\delta_s^{-2}C(\delta_s, \Sigma
^{p_s}_S(\lambda_s, \alpha_s, d_s), {\llVert \cdot\rrVert} )\}/\{\delta
_t^{-2}C(\delta_t, \Sigma^{p_t}(L_t, \break \alpha_t, d_t), \llVert \cdot
\rrVert )\}
$ and $B = 1 + \max_s d_s(\bar d - 1)$. 
\end{theorem}

\begin{pf}
Fix $k$ mutually exclusive subsets $B_1, \ldots, B_k$ of $\{1,\ldots,p\}
$ with $\llvert B_s\rrvert = p_s$, $1 \le s \le p$. Let $\Sigma^s$
denote the
space of norm $\lambda_s$, $\alpha_s$-smooth regression functions
that select $d_s$ predictors from $B_s$ and none from the other
subsets, that is, $\Sigma^s = \bigcup_{b \in\{0,1\}^p, \llvert
b\rrvert = d_s, \operatorname{support}(b)
\subset B_s}T^b(\lambda_s C^{\alpha_s, d_s}_1)$. These
subsets are mutually orthogonal since $f \in\Sigma^s$ and $f' \in
\Sigma^t$, $s \ne t$ pick disjoint sets of predictors and $f, f' \in
\mathcal{Z}_p$. Clearly, $\bigoplus_{s = 1}^k \Sigma^s \subset\Sigma
^{p,k,\bar d}_A(\lambda, \alpha, d)$. Let $f_i = \sum_{s = 1}^k
f_{is}$, $i = 1, \ldots, N$, be a \mbox{$\llVert \delta\rrVert /2$-}packing
set of
$\bigoplus_{s = 1}^k \Sigma^s$ under ${\llVert \cdot\rrVert} $. We
must have
%
\begin{equation}
N \ge\frac{1}4 \Biggl\{\frac{3}{4C^*} \sum
_{s = 1}^k C\bigl(\delta_s,
\Sigma^{p_s}_S(\lambda_s,
\alpha_s, d_s), \llVert\cdot\rrVert\bigr) - k \log2
\Biggr\},
\end{equation}
by an application of Lemma \ref{lemadd} with $c = 1/2$, coupled with
the fact that $\Sigma^s$ is isomorphic with $\Sigma^{p_s}_S(\lambda
_s, \alpha_s, d_s)$. Also, by Lemma \ref{lembasic} and the packing
set construction used in the proof of Lemma \ref{lemadd}, each
$f_{si}$ can be chosen to belong to $\Sigma^s \cup C^\infty(\bbR
^p)$. Define $g_1, \ldots, g_N$ as: $g_i(x) = \Delta^{\bar\alpha}
f_i(x / \Delta)$ where $\bar\alpha= \max_s \alpha_s$. Then each
$g_i \in\Sigma^{p,k,\bar d}_A(\lambda, \alpha, d)$ and, $\llVert g_i -
g_j\rrVert _Q \ge\underline q^{1/2} \Delta^{\bar\alpha+ \max_s
d_s/2}\llVert
f_i - f_j\rrVert $, since every $f_{is} - f_{js}$ involve at most $\max_s
d_s$ many variables and they are orthogonal across $s$. This proves the
first assertion of the theorem.

In light of the well-known relation $V(\varepsilon, A, \llVert \cdot
\rrVert ) \le
C(\varepsilon, A, {\llVert \cdot\rrVert} ) \le V(\varepsilon/2,\break  A,
{\llVert \cdot\rrVert} )$ between
packing and covering entropies of subsets in a metric space, and the
fact that ${\llVert \cdot\rrVert} _Q \le\bar q^{1/2} \llVert \cdot
\rrVert $, the second
assertion can be established by showing
\[
V\bigl(2\sqrt{B}\llVert\delta\rrVert, \Sigma_A^{p, k, \bar d}(
\lambda, \alpha, d), \llVert\cdot\rrVert\bigr) \le\sum
_{s = 1}^k V\bigl(\delta_s,
\Sigma^p_S(\lambda_s,
\alpha_s, d_s), \llVert\cdot\rrVert\bigr).
\]
For every $1 \le s \le k$, let $\mathcal{C}^s$ be a minimal $\delta
_s/\lambda_s$-covering set of $C^{\alpha_s, d_s}_1$. For each $s$,
replace every element $f \in\mathcal{C}^s$ by its centered version
$\bar f = f - \int f(x) \,dx$. The new $\mathcal{C}^s$ remains a
$2\delta_s/\lambda_s$-covering set of $C^{\alpha_s, d_s}_1 \cap
\mathcal{Z}_{d_s}$. Take
\[
\mathcal{C}_A = \Biggl\{f = \sum_{s = 1}^k
\lambda_s T^{b^s} f_s\dvtx f_s \in
\mathcal{C}^s, 1 \le s \le k, \bigl(b^1, \ldots,
b^k\bigr) \in\mathcal{B}^{p,k,d,\bar d} \Biggr\}.
\]
Any $f \in\Sigma_A^{p,k,\bar d}(\lambda, \alpha, d)$ equals $f =
\sum_s \lambda_s T^{b^s} f_s$ for some $f_s \in C^{\alpha_s,d_s}_1
\cap\mathcal{Z}_{d_s}$, $1 \le s \le k$ and $(b^1,\ldots,b^k) \in
\mathcal{B}^{p,k,d,\bar d}$. Find $f^*_s \in\mathcal{C}^s$ such that
$\llVert f_s - f^*_s\rrVert \le2\delta_s/\lambda_s$, $1 \le s \le k$
and set
$f^* = \sum_s \lambda_s T^{b^s} f^*_s \in\mathcal{C}_A$. Since
every $f_s - f^*_s \in\mathcal{Z}_{d_s}$, we get $\langle T^{b^s}(f_s
- f^*_s), T^{b^t}(f_t - f^*_t)\rangle= 0$ whenever $\sum_{j = 1}^p
b^s_j b^t_j = 0$, that is, $b^s$, $b^t$ have no shared selection. By
assumption on $\mathcal{B}^{p,k,d,\bar d}$, for every $s$, there are
at most $d_s(\bar d - 1)$ many $t \ne s$ with shared selection.
Therefore, by Lemma \ref{lemoverlap}, $\llVert f - f^*\rrVert ^2 \le
B \sum_{s =
1}^k \lambda_s^2 \llVert f_s - f^*_s\rrVert ^2 \le4B \llVert \delta
\rrVert ^2$.
Consequently, $\mathcal{C}_A$ gives a $(2\sqrt B \llVert \delta\rrVert
)$-covering of $\Sigma^{p,k,\bar d}_A(\lambda, \alpha, d)$. This\vspace*{1pt}
completes the proof since $V(\delta_s, \Sigma^p_S(\lambda_s, \alpha
_s, d_s), {\llVert \cdot\rrVert} ) \ge\log\llvert \mathcal
{C}^s\rrvert $ for every $1 \le s \le k$.
\end{pf}

%
%
\begin{lemma}
\label{lemoverlap}
Suppose $f_1, \ldots, f_k$ are elements of a Hilbert space $(\bbH,
{\llVert \cdot\rrVert} _{\bbH})$ and for any $1 \le s \le k$, let $r_s
= \llvert \{1 \le t
\le k\dvtx \langle f_s, f_t\rangle_\bbH\ne0\}\rrvert $. Then
$\llVert \sum_{s = 1}^k
f_s\rrVert _\bbH^2 \le\max_s r_s\sum_{s = 1}^k \llVert f_s\rrVert
^2_\bbH$.
\end{lemma}

\begin{pf}
Since $2\langle f, g \rangle_\bbH\le\llVert f\rrVert ^2_\bbH+
\llVert g\rrVert ^2_\bbH$,
we have
\begin{eqnarray*}
\biggl\llVert{ \sum_s} f_s\biggr
\rrVert^2_\bbH &=& \mathop{\sum}_{s,t}
\langle f_s, f_t \rangle_\bbH
\le \frac{1}2 \mathop{\sum}_{\langle f_s,
f_t\rangle_\bbH\ne0}\bigl(\llVert
f_s\rrVert^2_\bbH+ \llVert f_t
\rrVert^2_\bbH\bigr)
\\
&\le& \max_s
r_s \sum_s \llVert f_s
\rrVert^2_\bbH.
\end{eqnarray*}\upqed
\end{pf}

%
%
\begin{lemma}\label{lenormc}
Suppose $\mathcal{F} \subset C(\bbR^p)$ satisfies $\sup_{f\in
\mathcal{F}}\llVert f\rrVert _{\infty}\leq1$. 
Then, for any sequence $\delta_n$ satisfying $n\delta_n^2\geq2
N(\delta_n,\mathcal{F},{\llVert \cdot\rrVert} _{\infty})$ and $\sum
_{n=1}^{\infty} e^{-n\delta_n^2}< \infty$,
\[
Q^\infty\Bigl( \Bigl\{x^{1\dvtx\infty}\dvtx \sup_{f\in\mathcal{F},
\llVert f\rrVert
_Q\leq\delta_n}
\llVert f\rrVert_n\geq4\delta_n \mbox{ infinitely often}
\Bigr\} \Bigr)= 0.
\]
\end{lemma}

\begin{pf}
Take $f\in\mathcal{F}$ and suppose that $\llVert f\rrVert _Q\leq\delta
_n$ and
$\llVert f\rrVert _n\geq4\delta_n$. Let $\{f_1,\ldots,f_N\}$ form an minimal
$\delta_n$-covering of $\mathcal{F}$ under the sup-norm with $2\log
N\leq n\delta_n^2$. Then there exists some $j_0\in\{1,\ldots,N\}$
such that $\llVert f-f_{j_0}\rrVert _{\infty}\leq\delta_n$. By the assumptions
on $f$, we have $\llVert f_{j_0}\rrVert _n\geq3\delta_n$ and
$\llVert f_{j_0}\rrVert _Q\leq
2\delta_n$, implying $\llvert \llVert f_{j_0}\rrVert _n^2-\llVert
f_{j_0}\rrVert _{Q}^2
\rrvert \geq5\delta_n^2$. By Bernstein's inequality, we have
\[
P \bigl\{\bigl\llvert\llVert f_{j_0}\rrVert_n^2-
\llVert f_{j_0}\rrVert_{Q}^2\bigr\rrvert\geq5
\delta_n^2 \bigr\}\leq2\exp\biggl[-\frac{5}{8}n
\delta_n^2 \biggr].
\]
Since there are at most $N$ choices for $j_0$, we get
\begin{eqnarray*}
P \Bigl\{\sup_{f\in\mathcal{F}, \llVert f\rrVert _Q \leq \delta
_n}\llVert f\rrVert_n\geq4
\delta_n \Bigr\} &\leq& \sum_{j=1}^NP
\bigl\{\bigl\llvert\llVert f_{j}\rrVert_n^2-
\llVert f_{j}\rrVert_{Q}^2\bigr\rrvert\geq5
\delta_n^2 \bigr\}
\\
&\leq& 2N\exp\biggl[-\frac{5}{8}n\delta_n^2
\biggr]\leq2\exp\biggl[-\frac
{1}{8}n\delta_n^2
\biggr],
\end{eqnarray*}
from which the results follows by the Borel--Cantelli lemma.
\end{pf}

\section*{Acknowledgments} The authors would like to thank the
Associate Editor and two referees for their insightful comments and
suggestions.


%

\printaddresses
\end{document}